\documentclass{amsart}
\usepackage{amsmath}
\usepackage{amssymb}
\usepackage{amscd}
\usepackage{amsthm}
\usepackage[all]{xy}
\usepackage[mathcal]{euscript}

\DeclareMathOperator{\Irr}{Irr}
\DeclareMathOperator{\Hom}{Hom}
\DeclareMathOperator{\spn}{span}

\newcommand{\cF}{\mathcal{F}}
\newcommand{\cG}{\mathcal{G}}
\newcommand{\cM}{\mathcal{M}}

\newcommand{\cT}{\mathcal{T}}
\newcommand{\cU}{\mathcal{U}}
\newcommand{\C}{\mathbb{C}}
\newcommand{\F}{\mathbb{F}}
\newcommand{\Z}{\mathbb{Z}}
\newcommand{\Ab}{\bar{A}}
\newcommand{\Hb}{\bar{H}}
\newcommand{\Ib}{\bar{I}}
\newcommand{\ba}{\mathbf{a}}
\newcommand{\bb}{\mathbf{b}}
\newcommand{\bc}{{\mathbf{c}}}

\newcommand{\Orb}{\mathcal{O}}

\newtheorem{thm}{Theorem}[section]
\newtheorem{prop}[thm]{Proposition}
\newtheorem{lem}[thm]{Lemma}
\newtheorem{cor}[thm]{Corollary}

\theoremstyle{definition}
\newtheorem{defn}[thm]{Definition}
\newtheorem{rmk}[thm]{Remark}
\newtheorem{rmks}[thm]{Remarks}
\theoremstyle{remark}

\newcounter{partsctr}
\newenvironment{parts}{\setcounter{partsctr}{1}\catcode`@=11%
 \def\item{\ifnum\value{partsctr}>1\par\fi(\arabic{partsctr})\qquad%
 \addtocounter{partsctr}{1}\def\@currentlabel{\arabic{partsctr}}}}{}

\title[Special pieces, Springer correspondence, and unipotent characters]{On special pieces, the Springer correspondence, and unipotent characters}
 
\author{Pramod N.~Achar}
\address{Department of Mathematics\\
  Louisiana State University\\
  Baton Rouge, LA 70803}
\email{pramod@math.lsu.edu}
\thanks{The research of the first author was partially supported by NSF grant~ DMS-0500873.}
\author{Daniel S.~Sage}
\email{sage@math.lsu.edu}
\thanks{The research of the second author was partially supported by NSF grant~ DMS-0606300.}

\begin{document}

\begin{abstract}
  Let $G$ be a connected reductive algebraic group over the algebraic
  closure of a finite field $\F_q$ of good characteristic.  In this
  paper, we demonstrate a remarkable compatibility between the
  Springer correspondence for $G$ and the parametrization of unipotent
  characters of $G(\F_q)$.  In particular, we show that in a suitable
  sense, ``large'' portions of these two assignments in fact coincide.
  This extends earlier work of Lusztig on Springer representations
  within special pieces of the unipotent variety.
\end{abstract}
 
\maketitle

\section{Introduction}
\label{sect:intro}

Let $G$ be a connected reductive algebraic group over the algebraic
closure of a finite field $\F_q$ of good characteristic, and let $W$
be its Weyl group.  The Springer correspondence, which we denote by
$\nu$, assigns to each irreducible representation of $W$ an
irreducible equivariant local system on a unipotent class of $G$.  On
the other hand, part of Lusztig's parametrization of unipotent
characters of finite reductive groups (or of his parametrization of
unipotent character sheaves) is a map $\kappa$ assigning to each
irreducible representation of $W$ an element of a certain finite set,
whose members may be thought of as irreducible equivariant local
systems on certain finite groups.

The main goal of this paper is demonstrate a certain compatibility between
these two maps.  Both maps assign to each representation of $W$ an orbit
in some space and a representation of an isotropy group.  Our main
result asserts that if we restrict our attention to a given two-sided
cell (in a sense to be made precise later), then there
is a natural order-preserving bijection between the two sets of orbits, and
that the relevant isotropy groups are canonically isomorphic. Furthermore,
we show that if we then identify corresponding orbits and isotropy groups,
then the Springer correspondence and Lusztig's map $\kappa$ actually
coincide.

This theorem generalizes various results due to Lusztig~\cite{Lnotes}.
For exceptional groups, its proof amounts to poring over tables
from~\cite{Ca}.  We will indicate in Section~\ref{sect:back} exactly
what poring is to be done, but we will not tabulate the results here.
(Some of the calculations required for the exceptional groups case
were first carried out by Kottwitz.)  The proof for classical groups is
carried out in Sections~\ref{sect:outline} and~\ref{sect:combin}.  In
Section~\ref{sect:outline}, we discuss various general phenomena in
the classical groups and produce an outline of a proof that is
independent of type, but of course lacking in combinatorial details.
Section~\ref{sect:combin} supplies these details in full for type $C$
and somewhat more sparsely in types $B$ and $D$.

It has often been the case that `combinatorial coincidences' in
representation theory hint at deep underlying relationships between
superficially disparate phenomena, and one wonders whether this is the
case here.  Since the Springer correspondence is intimately related to
the theory of character sheaves, which enjoy a classification similar
to that of characters of finite reductive groups, perhaps there is a
purely geometric proof of the results of this paper in the character
sheaf setting.  (A small hint about the geometry in the problem
appears in Lusztig's paper~\cite{Lnotes}: he comments on the
relationship between the main results of that paper and a theorem of
Kraft--Procesi~\cite{KP} on the geometry of special pieces in the
classical groups.) Such a proof would certainly shed light on why the
results of this paper are true; perhaps it would also shed light on
other aspects of the relationship among character sheaves, unipotent
characters, and two-sided cells.

We are grateful to Robert Kottwitz for several valuable discussions as
well as for sharing his exceptional groups calculations with us.

\section{Background and definitions}
\label{sect:back}

\subsection{Unipotent classes and Lusztig's canonical quotient group}

Let $\cU$ denote the unipotent variety in $G$.  For any $G$-stable
subvariety $Z$ of $\cU$, let $\Orb(Z)$ be the set of unipotent classes
contained in $Z$.  For each class $C \in \Orb(\cU)$, we denote by $A(C)$
the component group of the centralizer of some point $u \in C$, that is,
$A(C) = G^u/(G^u)^\circ$.  (Different choices of $u$ clearly result in
isomorphic groups $A(C)$; moreover, the isomorphism is determined up to an
inner automorphism.)  With $Z$ as above, let
\[
X(Z) = \{(C,\rho) \mid
\text{$C \in \Orb(Z)$ and $\rho \in \Irr(A(C))$}\}.
\]
Of course, this set can be identified with the set of irreducible
$G$-equivariant local systems on $G$-orbits in $Z$. 

In particular, we may regard the Springer correspondence as a map $\nu:
\Irr(W) \to X(\cU)$.  If $(C,\rho)$ is in the image of $\nu$, we write
$\chi_{C,\rho}$ for its preimage. 

Next, for a given unipotent class $C$, let
\[
K = \bigcap_{\substack{\rho \in \Irr(A(C)) \\
\text{$\chi_{C,\rho}$ and $\chi_{C,1}$ are in the} \\
\text{same two-sided cell}}} \ker \rho
\qquad\text{and}\qquad
\Ab(C) = A(C)/K.
\]
This group is called \emph{Lusztig's canonical quotient} of $A(C)$.
We further define
\[
X^0(Z) = \{(C,\rho) \mid
\text{$C \in \Orb(Z)$ and $\rho \in \Irr(\Ab(C))$}\}.
\]
This set can be identified with a subset of $X(Z)$, by pulling back
representations of $\Ab(C)$ to $A(C)$. 

\subsection{Special pieces}
\label{sect:specialpiece}

A \emph{special piece} is, by definition, the union of a special unipotent
class $C_0$ and those unipotent classes $C$ in its closure with the
property that for any other special class $C_1 \subsetneq \overline C_0$,
we have $C \not\subset \overline{C_1}$.  By a result of
Spaltenstein~\cite{Sp}, the full unipotent variety is the disjoint union of
the various special pieces.  Evidently, there is a one-to-one
correspondence between special pieces and two-sided cells of $W$: given a
two-sided cell $\bc$, we denote by $P_\bc$ the special piece determined by
the special unipotent class $C_0$ with the property that $\chi_{C_0,1}$ is
the unique special character in $\bc$. 

Let $\Irr(W)_\bc$ denote the set of irreducible characters of $W$ belonging
to $\bc$.  (This set is a ``family'' of
characters as defined in~\cite[\S 4.2]{Lbook}.)  Next, let
\[
\Irr(W)^0_\bc = \{ \chi \in \Irr(W)_\bc \mid
\text{$\nu(\chi) = (C,\rho)$ where $C \subset P_\bc$} \}.
\]
In other words, $\Irr(W)^0_\bc$ contains those characters in $\bc$
that are sent by $\nu$ to the ``correct'' special piece.  We also let
$\Irr(W)^0_{\bc,C}$ denote the set of characters in $\Irr(W)^0_\bc$
corresponding to local systems on $C$.  Lusztig has shown that a given
unipotent class $C$ is contained in $P_\bc$ if and only if $\chi_{C,1}
\in \Irr(W)_\bc$.\cite{Lnotes} Thus, if $\chi = \chi_{C,\rho} \in
\Irr(W)^0_\bc$, then $\chi_{C,\rho}$ and $\chi_{C,1}$ both belong to
$\bc$, so $\rho$ descends to a representation of $\Ab(C)$ by its
definition as a quotient of $A(C)$.  We conclude that if $\chi \in
\Irr(W)^0_\bc$, then $\nu(\chi) \in X^0(P_\bc)$.  (We will see later
that the converse is also true, namely that if $(C,\rho)\in
X^0(P_\bc)$ is in the image of the Springer correspondence, then
$\chi_{C,\rho}\in \Irr(W)^0_\bc$.)

\subsection{The parametrization of unipotent characters and character
sheaves}

Let $\bc$ be a two-sided cell of $W$.  In~\cite{Lbook}, Lusztig defined a
certain finite group $\cG_\bc$  attached to $\bc$.  (In the notation
of~\cite{Lbook}, this group is called $\cG_\cF$, where $\cF$ is a
family of characters.)  Let $Z
\subset \cG_\bc$ be a subset that is stable under conjugation.  By analogy
with the unipotent variety of $G$, we define $\Orb(Z)$ to be the set of
$\cG_\bc$-orbits (i.e., conjugacy classes) in $Z$, and for any $y \in
\Orb(\cG_\bc)$, we define $I(y)$ to be the centralizer $\cG_\bc^x \subset
\cG_\bc$ of some $x \in y$.  (Again, $I(y)$ is uniquely determined up to
inner automorphism).  We further define
\[
M(Z) = \{ (y, \sigma) \mid \text{$y \in \Orb(Z)$ and $\sigma \in
\Irr(I(y))$} \}.
\]
Regarding $Z$ as a discrete topological space, we may think of $M(Z)$ as
the set of isomorphism classes of irreducible $\cG_\bc$-equivariant
coherent sheaves on orbits in $Z$.

If $Z = \cG_\bc$, then the set $M(\cG_\bc)$ that we have described
coincides with the set $\cM(\cG_\cF)$ introduced by Lusztig
in~\cite{Lbook}.  He also constructed there an explicit
embedding of the set of irreducible representations in $\bc$ into
$\cM(\cG_\bc)$.  We denote this map by
\[
\kappa: \Irr(W)_\bc \hookrightarrow M(\cG_\bc).
\]

In the case where $G$ is split, the set of all unipotent representations
of $G(\F_q)$ is parametrized by the disjoint union of all $M(\cG_\bc)$ as
$\bc$ ranges over two-sided cells of $W$. On the other hand, those unipotent
representations arising by induction of the trivial character of a rational
Borel subgroup are parametrized by irreducible representations of $W$.  The
map $\kappa$ gives the relationship between these two parametrizations.

The same map $\kappa$ also arises in the parametrization of unipotent
character sheaves on $G$ (see~\cite[\S 17.5]{LCS4}).  Thus, the main
result of the paper may be interpreted in this setting as well.

\subsection{A topology on $\Ab(C)$ and $\cG_\bc$}
\label{sect:topl}

We now fix a unipotent class $C$.  In~\cite{Ad}, it was shown that the set of
conjugacy classes of $\Ab(C)$ admits a natural partial order (in which the
conjugacy class of the identity is the minimal element); moreover, $\Ab(C)$
has a canonical structure as a Coxeter group up to conjugacy and,
therefore, a well-defined class of ``parabolic subgroups.''  (It is well
known that the set of simple reflections of any irreducible finite Coxeter
group is unique up to conjugacy.  The point here is that $\Ab(C)$ is not, in
general, irreducible, but it nevertheless has a canonical decomposition
into irreducible factors.)  To any element $x \in \Ab(C)$, we
associate a parabolic subgroup $H_x$ via
\[
H_x = \text{the smallest parabolic subgroup of $\Ab(C)$ containing $x$}.
\]
(This definition makes sense because the set of parabolic subgroups is
closed under intersection, so $H_x$ is the intersection of all parabolic
subgroups containing $x$.)  Evidently, if $x_1$ and $x_2$ are conjugate,
then so are $H_{x_1}$ and $H_{x_2}$.

We endow $\Ab(C)$ with a topology by declaring a set to be open if it is a
union of parabolic subgroups.  In this topology, we see that for any two elements
$x_1, x_2 \in \Ab(C)$, the conjugacy class of $x_1$ is
contained in the closure of the conjugacy class of $x_2$ if and only if
$H_{x_2}$ is conjugate to a subgroup of $H_{x_1}$.  In particular, the
singleton consisting of the identity element is a dense open set.

Next, given a two-sided cell $\bc$, recall from~\cite{Lbook} that the
finite group $\cG_\bc$ can be identified with $\Ab(C_\bc)$, where $C_\bc$
is the unique special class in $P_\bc$.  We endow $\cG_\bc$ with a
topology via this identification, and we say that a $\cG_\bc$-equivariant sheaf $\cF$ on
$\cG_\bc$ is \emph{locally trivial} if, for each parabolic subgroup $H
\subset \cG_\bc$, the action of $H$ on $\cF(H)$ is trivial.  The stalk at
any point $x$ of a locally trivial $\cG_\bc$-equivariant sheaf then carries an
action of $\cG_\bc^x/(H_x \cap \cG_\bc^x)$.

For any conjugacy class $y \in \Orb(\cG_\bc)$, we define
\[
\Ib(y) = \cG_\bc^x/(H_x \cap \cG_\bc^x) \simeq N(H_x)/H_x
\qquad\text{for some $x \in y$;}
\]
as with $I(y)$, this group is determined independently of the choice of
$x$ up to inner automorphism.  Given a $\cG_\bc$-stable subset $Z
\subset \cG_\bc$, the set of isomorphism classes of irreducible locally
trivial $\cG_\bc$-equivariant sheaves on $Z$ is parametrized by the set
\[
M^0(Z) = \{ (y, \sigma) \mid \text{$y \in \Orb(Z)$ and $\sigma \in
\Irr(\Ib(y))$} \}.
\]

\subsection{Statement of the main result}

We regard the set $\Orb(P_\bc)$ as being partially ordered in the usual way: $C_1 \le C_2$ if $C_1 \subset \overline{C_2}$.  In addition, we also endow $\Orb(\cG_\bc)$ with an analogous partial order, using the topology defined above.

\begin{thm}\label{thm:reform}
Let $\bc$ be a two-sided cell in $W$.
\begin{enumerate}
\item There is an order-preserving injective map $t: \Orb(P_\bc) \to
  \Orb(\cG_\bc)$, characterized by the property that $(t(C),1) =
  \kappa(\chi_{C,1})$.  Moreover, the image of $t$ consists of the set
  of conjugacy classes of a subgroup $\cG'_\bc \subset
  \cG_\bc$.\label{pt:classes}
\item For each $C \in \Orb(P_\bc)$, there is a canonical isomorphism
  $\Ab(C) \simeq \Ib(t(C))$, and hence a canonical bijection $T:
  X^0(P_\bc) \to M^0(\cG'_\bc)$.\label{pt:group}
\item The following diagram commutes:\label{pt:diag}
\[
\xymatrix@R=10pt@C=5pt{
& *+{\raisebox{7pt}{$\Irr(W)^0_\bc$}}
\ar@{_{(}->}[dl]_{\nu}\ar@{^{(}->}[dr]^{\kappa} \\
X^0(P_\bc) \ar[rr]^{\sim}_{T} && M^0(\cG'_\bc)
}
\]
\item In each simple type except $G_2$, $F_4$, and $E_8$, the maps
  $\nu$ and $\kappa$ in the diagram are bijections for each two-sided
  cell.  For the three remaining types, there is exactly one cell for
  which the maps are not surjective, corresponding to the special
  unipotent classes $G_2(a_1)$, $F_4(a_3)$, and $E_8(a_7)$ (using the
  notation of \cite{Ca}). In each case, a single local system  is
  missing from the image of $\nu$, namely the special class
  together with the sign representation of $\Ab(C)$ (which is $S_3$,
  $S_4$, and $S_5$ respectively).\label{pt:bij}
\end{enumerate}
\end{thm}

It suffices to prove the theorem for each simple root system.
Moreover, there are no nontrivial local systems in type $A$, so the
theorem is trivial in this case.
 
\begin{proof}[Proof for the exceptional types]
  Parts~(\ref{pt:classes}) and~(\ref{pt:group}) of this theorem were
  already proved in Theorem~0.4 and Proposition~0.7 of~\cite{Lnotes}.
  (Proposition~0.7(b) of {\it loc. cit.} asserts that the group we
  have called $\Ib(t(C))$ (there, $N(H_C)/H_C$) is isomorphic to
  $A(C)$, not $\Ab(C)$.  However, that statement is only meant to
  apply to \emph{nonspecial} classes in the exceptional groups.  By
  comparing the tables in~\cite{Som} with those in~\cite{Ca}, one
  verifies that for all such classes, it is in fact the case that
  $A(C) = \Ab(C)$.)

Only the commutativity of the diagram in part~(\ref{pt:diag}) and the
statements about the images of the maps $\nu$ and $\kappa$ in
part~(\ref{pt:bij}) remain.  For
this, one undertakes the laborious task of comparing tables for the
Springer correspondence in, say, \cite{Ca} with Lusztig's enumeration
of values of $\kappa$ in~\cite{Lbook}.
\end{proof}

\begin{rmks} \begin{enumerate} \item The commutative diagram in part~(\ref{pt:diag}) of the
  theorem restricts to give a diagram for each unipotent orbit
  $C\subset P_\bc$:
\[
\xymatrix@R=10pt@C=5pt{
& *+{\raisebox{7pt}{$\Irr(W)^0_{\bc,C}$}}
\ar@{_{(}->}[dl]_{\nu}\ar@{^{(}->}[dr]^{\kappa} \\
X^0(C) \ar[rr]^{\sim}_{T} && M^0(t(C))
}
\]
Again, the maps $\nu$ and $\kappa$ are bijections except in the three
cases mentioned above.
\item It is a corollary of the theorem and its proof that for any
  $\chi\in\Irr(W)$, $\chi \in \Irr(W)^0_\bc$ if and only if $\nu(\chi)
  \in X^0(P_\bc)$.  Indeed, this follows from part~(\ref{pt:bij}) of the
  theorem together with the injectivity of the Springer correspondence
  for all cells besides the three anomalies.  However, in these
  three cases, the single element in $X^0(P_\bc)$ not in
  $\nu(\Irr(W)^0_\bc)$ is not in the image of the Springer
  correspondence.
\end{enumerate}
\end{rmks} 

\section{The classical types: Outline}
\label{sect:outline}

The proof of the main theorem for the classical types $B$, $C$, and
$D$ comes down to calculations with various combinatorial objects that
parametrize $\Irr(W)$, $X(\cU)$, and $X(\cG_\bc)$.  The details differ
in the three types, but they can all be placed in a common framework,
which we will now develop.

\subsection{$\Ab(C)$ for classical groups}
\label{sect:Vsp}

Recall that for the classical root systems, the groups $\cG_\bc$ are
always products of copies of $\Z/2\Z$.  In particular, they are always
abelian.  Thus, the set $\Irr(\cG_\bc)$ can be identified with the
group $\hat\cG_\bc = \Hom(\cG_\bc, \C^\times)$, and $M(\cG_\bc)$ can,
in turn, be identified with $\cG_\bc \times \hat \cG_\bc$.  In this
section, we will give a particular realization of these spaces that is
well-suited to the combinatorics that will follow.  Furthermore,
we describe the sets $\Irr(\cG_\bc/H_x)$ with respect to this
realization.

Suppose $\cG_\bc \simeq (\Z/2\Z)^f$.  Let $\tilde V_\bc$ be the 
$(2f+1)$-dimensional $\F_2$-vector space with basis
\[
e_0, e_1, \ldots, e_{2f}
\]
and $V_\bc$ its quotient 
by the relation
\begin{equation}\label{eqn:Vrel}
e_0 + e_2 + e_4 + \cdots + e_{2f} = 0
\end{equation}
(so $\dim V_\bc = 2f$).
We endow $\tilde V_\bc$ with a symmetric bilinear form given by
\[
\langle e_i, e_j\rangle =
\begin{cases}
1 & \text{if $|i - j| = 1$,} \\
0 & \text{otherwise.}
\end{cases}
\]
This form has kernel spanned by $e_0 + e_2 + \cdots + e_{2f}$, so it
induces a nondegenerate symmetric bilinear form on $V_\bc$.  We make
the identification
\[
\cG_\bc = \spn\{e_1, e_3, \ldots, e_{2f-1}\},
\]
and we remark, in particular, that the canonical Coxeter generators referred to in Section~\ref{sect:topl} are precisely the elements $e_1, e_3, \ldots, e_{2f-1}$. Via the bilinear form $\langle,\rangle$, we can also identify
\[
\hat \cG_\bc = \Hom(\cG_\bc, \C^\times) = \spn \{e_0, e_2, \ldots, e_{2f}\},
\]
viewed as a subspace of $V_\bc$.  Note that every character has
exactly two representatives in $\tilde V_\bc$.

Now, take $x \in \cG_\bc$.  There is a set of integers $n(x) \subset \{1,3\ldots, 2f-1\}$ such that
\[
x = \sum_{j \in n(x)} e_j.
\]
Evidently, the parabolic subgroup $H_x$ is simply the span of $\{e_i \mid i\in
n(x)\}$.  To describe $\Irr(\cG_\bc/H_x)$ as a subset of $\hat \cG_\bc$, we
need a slightly refined description of $n(x)$.  An element $\lambda \in
\spn \{e_0, e_2, \ldots, e_{2f}\}$ is called a \emph{block} of $x$ if it is of the form
\[
\lambda = e_i + e_{i+2} + \cdots + e_{i+2k}
\]
where
\[
\{i+1,i+3,\ldots,i+2k-1\} \subset n(x), \quad
i-1 \notin n(x), \quad i+2k+1 \notin n(x).
\]
Note that $i$ must be even here.  It is permitted that $k = 0$, so a
single term $e_i$ is a block if $i$ is even and $i-1, i+1 \notin
n(x)$.  We will also view $\lambda$ as a character in $\hat \cG_\bc$
by taking its image in $V_\bc$.  Next, the element $m \in \cG_\bc$
given by
\[
m = e_{i+1} + e_{i+3} + \cdots + e_{i+2k-1}
\]
is called the \emph{component} of $x$ associated to the block $\lambda$. 
(If $x$ has a block that is a single term $e_i$, then the associated
component is $0$.)

Clearly, $x$ is the sum of its components, with distinct components
orthogonal and each $e_j$ for $j$ odd appearing in at most one
component.  Similarly, every $e_j$ with $j$ even appears in exactly
one block of $x$.  Thus, the sum of the blocks is $0$
(by equation~\eqref{eqn:Vrel}), and any proper subset of the blocks is linearly
independent.

\begin{lem}\label{lem:Ib}
  $\Irr(\cG_\bc/H_x)$ is the subspace of $\hat \cG_\bc$ spanned by the
  images of the blocks of $x$.
\end{lem}
\begin{proof}
Let $\lambda$ be a block of $x$, and let $m$ be its associated component. 
It is clear that $\langle\lambda, e_j\rangle = 0$ for each $e_j$ appearing
in $m$.  Moreover, neither $i-1$ nor $i+2k+1$ belongs to $n(x)$, so in
fact, $\langle \lambda, e_j \rangle = 0$ for all $j \in n(x)$.  Thus,
$\lambda|_{H_x}$ is indeed trivial. 

Conversely, let $\mu \in \hat\cG_\bc$ be such that some, but not all,
of the $e_j$ appearing in $\lambda$ appear in $\mu$.  (This is
independent of the choice of representative of $\mu$ in $\tilde
V_\bc$.) In other words, $\mu$ is not in the span of the blocks of
$x$.  This means that we can find a pair of consecutive terms in $\lambda$, say
$e_{2p}$ and $e_{2p+2}$, such that exactly one of them appears in
$\mu$.  Then $2p+1 \in n(x)$, yet it is clear that $\langle \mu,
e_{2p+1}\rangle = 1$.  Thus, $\mu|_{H_x} \ne 0$.
\end{proof}

\subsection{Symbols and $u$-symbols}

The combinatorial objects with which all our calculations in the
classical types will be carried out are called \emph{symbols} and
\emph{$u$-symbols}.  Both of these are certain finite sequences
$(a_0,a_1, \ldots, a_r)$ of nonnegative integers, satisfying various
additional conditions (that are different for symbols and for
$u$-symbols).  The value of $r$ ({\it i.e.}, the length of the
sequences), as well as the precise form of the additional conditions,
depend on the type.  The set of symbols and the set of $u$-symbols
both parameterize $\Irr(W)$ explicitly.  Accordingly, there is a
bijection
\[
i: \{\text{$u$-symbols}\} \to \{\text{symbols}\},
\]
which maps the $u$-symbol of an irreducible
representation to its symbol. 

\begin{rmk}
Our usage of the terms ``symbol'' and ``$u$-symbol'' is close to that
of~\cite{Lnotes}, but narrower than the usage elsewhere in the literature. 
The definitions that we will give correspond to what other authors refer to
as symbols or $u$-symbols ``of defect $1$'' (for types $B$ and $C$) or ``of
defect $0$'' (for type $D$). 

Also, symbols and $u$-symbols are usually written as arrays consisting
of two rows of numbers.  Here, following~\cite{Lnotes}, we will write them
simply as finite sequences, in which the even-numbered entries $a_0, a_2,
\ldots$ may be thought of as the ``upper row,'' and the odd-numbered
entries as the ``lower row.'' In all three types, symbols will have
strictly increasing upper and lower rows while $u$-symbols will have
$a_{i+2}-a_i\ge 2$.
\end{rmk}

Symbols and $u$-symbols both fall into families determined by their
entries. 

\begin{defn}
Two symbols $\ba$ and $\bb$ are said to be \emph{congruent} if they contain
the same entries (with the same multiplicities), possibly in different
orders.  A symbol $(a_0, \ldots, a_r)$ is \emph{special} if $a_0 \le \cdots
\le a_r$. 

Two $u$-symbols $\ba$ and $\bb$ are said to be \emph{similar} if they
contain the same entries (with the same multiplicities), possibly in
different orders.  A $u$-symbol $(a_0, \ldots, a_r)$ is
\emph{distinguished} if $a_0 \le \cdots \le a_r$. 
\end{defn}

In particular, each congruence class of symbols contains a unique special
symbol~\cite[Lemmas~1.6, 2.6, 3.6]{Lnotes}, and each similarity class of
$u$-symbols contains a unique distinguished symbol~\cite[\S 11.5]{LIC}.  

Many of the combinatorial results we require will be easier to state
and prove if we work with the indices of entries in a symbol or
$u$-symbol rather than the entries themselves.  In other words, we
make statements in terms of subsets of $\{0,\ldots, r\}$ rather than
subsequences of $(a_0, \ldots, a_r)$.  The following notation will be
particularly useful:
\[
[k,l] = \{k, k+1, \ldots, l\}.
\]

\begin{defn}
Let $\ba = (a_0, \ldots, a_r)$ be a $u$-symbol, and let $\ba' = i(\ba) =
(a'_0, \ldots, a'_r)$ be the corresponding symbol. 

A number $k \in [0,r]$ is called an \emph{isolated point} (of $\ba$ or
$\ba'$) if $a'_k \ne a'_l$ for all $l \ne k$. 

A subset $[k,l] \subset [0,r]$ is called a \emph{ladder} of $\ba$ if, for
some $a$, we have 
\begin{gather*}
(a_k, a_{k+1}, \ldots, a_l) = (a, a+1, \ldots, a+l-k), \\
\text{$a_{k-1} < a-1$ if $k > 0$,\qquad and\qquad $a_{l+1} > a+l-k+1$ if $l
< r$.} 
\end{gather*}
(Ladders are called intervals in some earlier references.)

A subset $[k,l] \subset [0,r]$ is called a \emph{staircase} of $\ba$ if,
for some $a$, we have 
\begin{gather*}
(a_k, a_{k+1}, \ldots, a_l) = (a, a, a+2, a+2,\ldots, a+l-k-1, a+l-k-1), \\
\text{$a_{k-2} < a-2$ if $k > 1$,\qquad and\qquad $a_{l+2} > a+l-k+1$ if $l
< r-1$.} 
\end{gather*}
(Note that $l-k+1$ is necessarily even, if $[k,l]$ is a staircase.)

A set $[k,l]$ that is either a ladder or a staircase is called a
\emph{part}.  The number $k$ is called the \emph{bottom} of the part,
and $l$ is called the \emph{top}.  The \emph{length} of the part is
defined to be $l-k+1$.
\end{defn}

We will only be interested in parts for distinguished $u$-symbols.  It
is easy to see that the set of indices for any finite nondecreasing
sequence of integers satisfying $a_{i+2}-a_i\ge 2$ is the disjoint
union of ladders and staircases.  In particular, for a distinguished
$u$-symbol, the full set of indices $[0,r]$ is the disjoint union of
all the parts.

We now describe a construction that will have multiple uses in the sequel.

\begin{defn}
Let $\ba = (a_0, \ldots, a_r)$ be a symbol (resp.~$u$-symbol), and let
$\mu$ be a subset of $\{0, \ldots, r\}$, such that for each $i \in \mu$,
the entry $a_i$ is distinct from all other entries of $\ba$.  Let $A =
\{a_0, a_2, a_4, \ldots\}$, $B = \{a_1, a_3, a_5, \ldots\}$, and $Z =
\{a_i \mid i \in \mu\}$.  Define two new sequences as follows:
\begin{align*}
A' = (b_0, b_2, b_4 \ldots) &=
\text{$(A \smallsetminus Z) \cup (B \cap Z)$, arranged in increasing
order,} \\
B' = (b_1, b_3, b_5 \ldots) &=
\text{$(B \smallsetminus Z) \cup (A \cap Z)$, arranged in increasing
order.}
\end{align*}
If the lengths of $A'$ and $B'$ are the same as the lengths of $A$ and
$B$, respectively, then we can combine them into a single sequence
$\bb = (b_0, b_1, b_2, \ldots, b_r)$, which is readily seen (for each
type) to be a new symbol (resp.~$u$-symbol).  We say that $\bb$ is
obtained by \emph{twisting} $\ba$ by $\mu$, and we denote it by $\bb =
\ba^\mu$.
\end{defn}

It is immediate that the twisting operation preserves congruence
classes of symbols and similarity classes of $u$-symbols.

\subsection{Symbols, cells, and $\kappa$}

Following~\cite{Lspec}, we know that two symbols are congruent if and only if
they correspond to representations of $W$ lying in the same two-sided cell.
 Moreover, the symbols we have called ``special'' are precisely those
corresponding to special representations of $W$. 

Now, note that any two congruent symbols have the same set of entries at
isolated points, so there is a natural bijection between their respective
sets of isolated points.  

\begin{defn}
Let $\ba_0$ be a special symbol, and let $\ba$ be a symbol in its
congruence class.  An isolated point of $\ba$ is said to be
\emph{displaced} if it has parity opposite to that of the corresponding
isolated point of $\ba_0$.  An isolated point of a $u$-symbol $\bb$ is
displaced if it is displaced in $i(\bb)$.  Every symbol (and
$u$-symbol) has an even number of displaced
isolated points. 
\end{defn}

Later, for each classical type, we will obtain a more direct
characterization of displaced isolated points. 

Let $\tilde V(\ba_0)$ be the power set of the set of isolated points,
regarded as an $\F_2$-vector space with addition given by symmetric
difference of sets.  Let $V(\ba_0)$ be the subspace of $\tilde V(\ba_0)$
consisting of sets with even cardinality.  We endow $V(\ba_0)$ with a
symmetric bilinear form by putting $\langle v, w \rangle = |v
\cap w| \mod 2$.  (Here $|v \cap w|$ is the cardinality of the set $v \cap
w$.) 

Since any symbol has an even number of displaced isolated points, 
there is a natural map
\[
\tilde\kappa:
\left\{\begin{array}{c}
\text{symbols}\\
\text{congruent to $\ba_0$}
\end{array}\right\}
\to V(\ba_0)
\qquad\text{via}\qquad
\ba \mapsto
\left\{\begin{array}{c}
\text{displaced}\\
\text{isolated points of $\ba$}
\end{array}\right\}
\]

Lusztig has shown that if $\bc$ is the two-sided cell corresponding to
the special symbol $\ba_0$, then there is a natural surjective map
$\pi$ from $V(\ba_0)$ to the space $V_\bc$ of Section~\ref{sect:Vsp};
in fact, it is an isomorphism (even an isometry) except in type
$D$.~\cite{Lbook} We will describe this map explicitly in each type by
labelling certain elements of $V(\ba_0)$ with the names $e_0, e_1,
\ldots, e_{2f}$.  With this in mind, the map $\tilde\kappa$ defined
above can be regarded as a combinatorial version of Lusztig's map
$\kappa: \Irr(W)_\bc \to V_\bc = M(\cG_\bc)$, in that the following
diagram commutes:
\begin{equation}\label{eqn:Ldiag}
\begin{array}{c}
\xymatrix{
*{\left\{\txt{symbols\\ congruent to $\ba_0$}\right\}}
\ar[rr]^{\tilde\kappa} \ar@{=}[d] && V(\ba_0) \ar[d]^{\pi} \\
\Irr(W)_\bc \ar[rr]^{\kappa} && V_\bc}
\end{array}
\end{equation}
(For a detailed discussion,
see~\cite[pp.~86--88]{Lbook} for types $B$ and $C$
and~\cite[pp.~92--94]{Lbook} for type $D$.)

It is not difficult to see that if $\ba$ is a symbol and $\mu$ is a set
consisting of an even number of isolated points, then
$\tilde\kappa(\ba^\mu) = \tilde\kappa(\ba) + \mu$ (assuming $\ba^\mu$ is
defined).  A more significant observation is that, by the commutativity of the
above diagram, we have 
\begin{equation}\label{eqn:Lcommute}
\kappa(\ba^\mu) = \kappa(\ba) + \pi(\mu).
\end{equation}

\subsection{$u$-Symbols and the Springer correspondence}

Following~\cite{LIC}, we know that two $u$-symbols are similar if and only
if the Springer correspondence maps the corresponding representations to
local systems on the same unipotent class.  Moreover, the distinguished
$u$-symbols are those corresponding to Springer representations (that is,
representations mapped to the trivial local system on some unipotent
class).  There is thus a bijection between distinguished $u$-symbols and
unipotent classes. 

Given a distinguished $u$-symbol $\ba = (a_0, \ldots, a_r)$, let
$H(\ba)$ be the power set of the set of ladders of $\ba$, regarded as
an $\F_2$-vector space with addition given by symmetric difference of
sets.

Let $C$ be the unipotent class corresponding to $\ba$.  In each classical
type, there is a surjective map $p: H(\ba) \to \Irr(A(C))$, such that the
set of all ladders lies in its kernel. (To be more precise, $\Irr(A(C))$
can be canonically identified with a quotient of some subspace of
$H(\ba)$, as described in~\cite[pp. 419--423]{Ca}.  The map $p$ may not be canonically
determined, but we may choose it so that its kernel contains the set of all
ladders.)

By abuse of language, we will often regard an element $\mu \in
H(\ba)$ as a subset of $[0,r]$ rather than as a set of ladders, by
replacing $\mu$ by the union of its members.  This abuse allows us to
speak of twisting a $u$-symbol by an element of $H(\ba)$.  

It turns out that all $u$-symbols in the similarity class of $\ba$
arise by twisting $\ba$ by a suitable set of ladders.\cite{LIC}
Moreover, this set of ladders is uniquely determined.  Given a
$u$-symbol $\ba'$ similar to $\ba$, we denote by $\tilde\nu(\ba') \in
H(\ba)$ the unique element such that $\ba' = \ba^{\tilde\nu(\ba')}$.
The following analogue of~\eqref{eqn:Ldiag} commutes, so $\tilde\nu$
can be thought of as a combinatorial version of the Springer
correspondence:
\begin{equation}\label{eqn:Sdiag}
\begin{array}{c}
\xymatrix{
*{\left\{\txt{$u$-symbols\\ similar to $\ba$}\right\}}
\ar[rr]^{\tilde\nu} \ar@{=}[d] && H(\ba) \ar[d]^{p} \\
\Irr(W)_C \ar[r]^-{\nu} & C\times\Irr(A(C))\ar[r]^-{\cong}&\Irr(A(C))}
\end{array}
\end{equation}
(A detailed account is given in~\cite[\S 12]{LIC} for type $C$
and~\cite[\S 13]{LIC} for types $B$ and $D$.)  Here, $\Irr(W)_C
\subset \Irr(W)$ is just the set of representations to which $\nu$
assigns local systems on $C$. In addition, we have an analogue
of~\eqref{eqn:Lcommute}:
\begin{equation}\label{eqn:Scommute}
\nu(\ba^\mu) = (C, p(\mu)),
\end{equation}
whenever $\ba^\mu$ is defined.

We will also introduce subsets of $[0,r]$ called the \emph{blocks} of
$\ba$, which are certain unions of ladders and staircases.  Each block
will contain at least one ladder, and each ladder will be contained in
a block.  Moreover, distinct blocks will have empty intersection.  Let
$B(\ba)$ be the set of blocks and $\Hb(\ba)$ the power set of
$B(\ba)$.  By assigning to a block the set of ladders contained within
it, we obtain an injective map from $\Hb(\ba)$ to $H(\ba)$, and we
will regard $\Hb(\ba)$ as a subset of $H(\ba)$ via this map.  In
particular, twisting by a set of blocks will mean twisting by the set
of ladders in these blocks.

Those $u$-symbols which may be obtained from $\ba$ upon twisting by a
set of blocks will be of particular interest.  Let $\cT(\ba)$ be the set
of those elements $\mu\in\Hb(\ba)$ for which $\ba^\mu$ is defined.  We
will see that it is a transversal of the two element subgroup
$\{\varnothing, B(\ba)\}$ and thus a set of size $2^{|B(\ba)|-1}$
(except in the trivial case where $B(\ba)$ is empty, which can only
occur in type $D$).

\subsection{Isolated points and blocks}

In this section, we give an outline of the proof of
Theorem~\ref{thm:reform} in the classical root systems.  Specifically, we
state a number of intermediate lemmas and propositions and indicate
how these statements together imply the theorem.  The proofs of most
of the intermediate results are slightly different in each type and
depend on the details of the combinatorial definitions of symbols and
$u$-symbols. Here, we only provide those proofs which are independent
of type.

We fix the following notations: let $\ba$ be a distinguished $u$-symbol,
and let $\ba_0$ be the unique $u$-symbol such that $i(\ba_0)$ is special
and congruent to $i(\ba)$.  (Note that $\ba_0$ is necessarily
distinguished, since all special representations of $W$ are Springer
representations.)  Let $\bc$ be the two-sided cell corresponding to
$i(\ba_0)$, and let $C$ be the unipotent class corresponding to $\ba$. 

The first result tells us where isolated points occur in a $u$-symbol:

\begin{prop}\label{prop:isol}
Let $\ba$ be a distinguished $u$-symbol. The top and bottom of each block
of $\ba$ are nondisplaced isolated points.  In addition to these, each
block contains an even number (possibly zero) of additional isolated
points, all of which are displaced.   
\end{prop}

We remark that the terms ``top'' and ``bottom'' as applied to blocks are
not quite the same as for ladders and staircases---the details vary by type
and will be given in Section~\ref{sect:combin}.  In particular, not every
block has both a top and a bottom, so this proposition does not imply that
there are necessarily an even number of isolated points.  Nevertheless, we
will see in each case that there is at most one block with an odd number of
isolated points.

Given a block $b$ of a distinguished $u$-symbol $\ba$, define $\lambda_b
\in V(i(\ba_0))$ by
\[
\lambda_b =
\begin{cases}
\text{the set of all isolated points in $b$} &
  \text{if that set has even cardinality;} \\
\text{the set of all isolated points \emph{not} in $b$} &
  \text{otherwise.}
\end{cases}
\]
The last sentence of the previous paragraph implies that $\lambda_b$ does
indeed belong to $V(i(\ba_0))$.  We extend the definition of
$\lambda_b$ to elements $b\in\Hb(\ba)$ by linearity.  Next, let $m_b \in V(i(\ba_0))$ be the set
of displaced isolated points of $b$.  (Accordingly, for most blocks, $m_b$ is
obtained from $\lambda_b$ by omitting the top and bottom of $b$.)

Note that according to the preceding proposition,
$\tilde\kappa(i(\ba)) = \sum m_b$, where $b$ runs over the blocks of
$\ba$.

\begin{prop}\label{prop:lusztig}
For any block $b$ of $\ba$, the element $\pi(m_b) \in \cG_\bc \times \hat
\cG_\bc$ is of the form $(x,1)$ for some $x \in \cG_\bc$.   
\end{prop}

As a consequence of this proposition, we recover Lusztig's result
\cite{Lnotes} that $\kappa(\chi_{C,1})$ is of the
form $(x,1)$ for some $x \in \cG_\bc$.  In particular, we henceforth
have available to us
the map $t: \Orb(P_\bc) \to \Orb(\cG_\bc)$ which takes the unipotent
orbit $C\subset P_\bc$ to the first coordinate of $\kappa(\chi_{C,1})$. 

\begin{prop}\label{prop:order}
If $C_1$ and $C_2$ are two unipotent classes in the same special piece,
then $t(C_1) \le t(C_2)$ in the partial order on $\cG_\bc$ if and only if 
$C_1$ is in the closure of $C_2$.
\end{prop}

Once this assertion is established, we will have completed the proof
of part~(\ref{pt:classes}) of Theorem~\ref{thm:reform}.  (The fact that
the image of $t$ is a subgroup of $\cG_\bc$ is proved in Theorem~0.4  of~\cite{Lnotes}.) 

\begin{prop}\label{prop:blocks}
  The elements $\pi(m_b)$, as $b$ runs over the blocks of $\ba$, are
  the components of $t(C)$, and the elements $\pi(\lambda_b)$ are the
  images in $V_\bc$ of the blocks of $t(C)$.
\end{prop}

Since the blocks of $\ba$ are a basis for $\Hb(\ba)$, and
$\Irr(\cG_\bc/H_x)$ is generated by the images of the blocks of $x$,
we see that the map $b \mapsto \lambda_b$ induces a surjective map
$\Hb(\ba) \to \Irr(\cG_\bc/H_{t(C)}) = \Irr(\Ib(t(C)))$.  Moreover, it
follows from the remark immediately preceding Lemma~\ref{lem:Ib} that
the kernel is generated by the sum of all the blocks of $\ba$.  We
will show that this is also the kernel of $p|_{\Hb(\ba)}: \Hb(\ba) \to
\Irr(A(C))$, but first, we must examine the set of $u$-symbols
obtainable from $\ba$ by twisting by a set of blocks.  We will see
that they correspond to the characters in $\Irr(W)^0_{\bc,C}$.  Recall
that $\cT(\ba)$ is the subset of $\Hb(\ba)$ consisting of those
elements $\mu$ for which $\ba^\mu$ is defined.

\begin{lem}\label{lem:transversal}
The set $\cT(\ba)$ contains at most one element from any coset of $\{\varnothing,B(\ba)\}$.
\end{lem}
\begin{proof} Suppose that both $\mu$ and $\mu^c=\mu+B(\ba)$ are in
  $\cT(\ba)$.  Since $B(\ba)$ is in the kernel of $p$,
  \eqref{eqn:Scommute} implies that $\nu(\ba^\mu)=\nu(\ba^{\mu^c})$,
  and the injectivity of the Springer correspondence shows that
  $\ba^\mu=\ba^{\mu^c}$.  This is a contradiction if $B(\ba)$ is
  nonempty, since a $u$-symbol $\ba'$ similar to $\ba$ is obtained by
  twisting $\ba$ by a unique set of ladders, namely $\tilde\nu(\ba')$.
  (The statement is trivial if $B(\ba)$ is empty.)
\end{proof}

\begin{prop}\label{prop:commute}
If $\ba$ is a distinguished $u$-symbol, then $i(\ba^b) =
(i(\ba))^{\lambda_b}$ for any $b\in \cT(\ba)$.  Furthermore, $\cT(\ba)$ is
a transversal in $\Hb(\ba)$ of the subgroup  $\{\varnothing,B(\ba)\}$.
\end{prop}
 
We now compute the kernel and range of the map $p|_{\Hb(\ba)}:
\Hb(\ba) \to \Irr(A(C))$.

\begin{prop}\label{prop:kernel}
The kernel of the restriction $p|_{\Hb(\ba)}: \Hb(\ba) \to \Irr(A(C))$ is
generated by $B(\ba)$, the sum of all the blocks of $\ba$.
\end{prop}
\begin{proof}
  Recall that $p: H(\ba) \to \Irr(A(C))$ was defined such that the set
  of all ladders is in its kernel.  Since every ladder is in a block,
  it is clear that the set of all blocks is in the kernel of
  $p|_{\Hb(\ba)}$.  On the other hand, distinct $\mu$'s in $\cT(\ba)$
  result in distinct twisted $u$-symbols, so the injectivity of the
  Springer correspondence and the commutativity of~\eqref{eqn:Sdiag}
  imply that $p|_{\Hb(\ba)}$ must be injective when restricted to
  $\cT(\ba)$.  But $\cT(\ba)$ is a transversal of
  $\{\varnothing,B(\ba)\}$, so the subgroup $\{\varnothing,B(\ba)\}$
  is indeed the kernel of $p|_{\Hb(\ba)}$.
\end{proof}

As a result, the map $p|_{\Hb(\ba)}: \Hb(\ba) \to \Irr(A(C))$ factors
through the surjection $\Hb(\ba) \to \Irr(\Ib(t(C)))$, inducing an
injective map $\Irr(\Ib(t(C))) \hookrightarrow \Irr(A(C))$.

Since $\Ab(C)$ is a quotient of $A(C)$, $\Irr(\Ab(C))$ is naturally a
subset of $\Irr(A(C))$.  In fact:

\begin{prop}\label{prop:imageAb}
The image of $p|_{\Hb(\ba)}: \Hb(\ba) \to \Irr(A(C))$ is precisely
$\Irr(\Ab(C))$.
\end{prop}

An immediate consequence of this last proposition is that the image of the
aforementioned injective map $\Irr(\Ib(t(C))) \hookrightarrow \Irr(A(C))$
is $\Irr(\Ab(C))$.  Part~(\ref{pt:group}) of Theorem~\ref{thm:reform} is
thus established; the map $T: X^0(P_\bc) \to M^0(\cG'_\bc)$ is described by
$T(C,p(b)) = (t(C), \pi(\lambda_b))$.

\begin{prop}\label{prop:bijection}
There is a bijection $\cT(\ba)\to\Irr(W)^0_{\bc,C}$ given by the map
$b\mapsto \ba^b$.  
\end{prop} 
\begin{proof} 
First, suppose that $b\in\cT(\ba)$, and consider the $u$-symbol
$\ba^b$.  It is similar to $\ba$ while Proposition~\ref{prop:commute}
shows that $i(\ba^b)$ is congruent to $i(\ba)$.  This says precisely
that $\ba^b$ corresponds to a character in $\Irr(W)^0_{\bc,C}$.

Now suppose that $\ba'$ is a $u$-symbol corresponding to a
representation in $\Irr(W)^0_{\bc,C}$.  Setting $\nu(\ba')=(C,\rho)$,
we have $\rho\in\Irr(\Ab(C))$, and since the image of $p|_{\Hb(\ba)}:
\Hb(\ba) \to \Irr(A(C))$ is precisely $\Irr(\Ab(C))$ (by
Proposition~\ref{prop:imageAb}), we can choose $b\in\bar H(\ba)$ with
$p(\mu)=\rho$.  Moreover, since $\cT(\ba)$ is a transversal for the
kernel of this map, we can choose $b\in\cT(\ba)$.  This means that
$\ba^b$ is defined, so by (5), $\nu(\ba^b)=(C,\rho)=\nu(\ba')$.  The
injectivity of $\nu$ gives $\ba'=\ba^b$.
\end{proof}

We obtain the following representation-theoretic corollary:

\begin{cor}\label{prop:Sprbijection} The Springer correspondence $\nu$ restricts to bijections
  $\Irr(W)^0_{\bc,C}\to X^0(C)$ and $\Irr(W)^0_\bc\to X^0(P_\bc)$.
\end{cor}
\begin{proof} 
  The fact that $\cT(\ba)$ is a transversal for the kernel of
  $p|_{\Hb(\ba)}$ implies that $|\cT(\ba)|=|\Irr(\Ab(C))|=|X^0(C)|$.
  The first bijection follows, and we get the second by taking
  the union over unipotent classes in $P_\bc$.
\end{proof}

We are now ready to complete the proof of the main theorem.  Given
$\chi=\chi_{C,\rho} \in \Irr(W)^0_\bc$, let $\ba'$ be its $u$-symbol
and $\ba$ the distinguished $u$-symbol for $C$.  By
Proposition~\ref{prop:bijection}, $\ba'=\ba^b$ for some
$b\in\cT(\ba)$.  We know that $\kappa(i(\ba)) = (t(C),1)$, so
equation~\eqref{eqn:Lcommute} and Proposition~\ref{prop:commute} give
$\kappa(i(\ba)^{\lambda_b}) = (t(C), \pi(\lambda_b))$.  In other words, we have $\kappa(\chi_{C,\rho}) = (t(C),
\pi(\lambda_b)) = T(C,\rho)$. Therefore,
\[
\kappa(\chi) = T(\nu(\chi))
\qquad
\text{for all $\chi \in \Irr(W)^0_\bc$.}
\]
This establishes the commutativity of the diagram in
part~(\ref{pt:diag}) of Theorem~\ref{thm:reform}.  Finally, we obtain
part ~(\ref{pt:bij}) of the theorem: the restriction of $\nu$ in the
diagram is a bijection by Proposition~\ref{prop:Sprbijection}, and so
$\kappa=T\circ\nu$ is a composition of bijections.  The proof of the
theorem is thus completed.

\section{The classical types: Combinatorics}
\label{sect:combin}

It remains, of course, to fill in a number of details separately in each of
the classical types.  The missing ingredients are as follows:
\begin{enumerate}
\item Definitions of ``symbol,'' ``$u$-symbol,'' and ``block.''
\item A careful study of where isolated points may occur.\label{pt:isol}
\item Proofs of Propositions~\ref{prop:isol}--\ref{prop:blocks},
  \ref{prop:commute}, and \ref{prop:imageAb}.
\end{enumerate}
Below, we give full details for each of these steps in type $C$.  In
particular, item~(\ref{pt:isol}) is treated in Lemmas~\ref{lem:isol}
and~\ref{lem:details}.  However, types $B$ and $D$ receive a much more
cursory treatment: we provide the requisite definitions, and we state
without proof the necessary results on isolated points (see
Lemmas~\ref{lem:isolB} and~\ref{lem:detailsB} for type $B$, and
Lemmas~\ref{lem:isolD} and~\ref{lem:detailsD} for type $D$), but we make
no comment at all on the other propositions, aside from giving an
explicit characterization (again without proof) of the transversal $\cT(\ba)$.  All
the missing proofs in types $B$ and $D$ are, of course, quite similar to
the corresponding proofs in type $C$.

\subsection{Type $C$}

Let $\Psi'_{2n,m}$ (resp.~$\Phi'_{2n,m}$) be the set of
sequences $\ba = (a_0, a_1, \ldots, a_{2m})$ satisfying the following
conditions:
\begin{align*}
0 &\le a_0, & 1 &\le a_1, & a_i &\le a_{i+2}-2, & \sum a_i &= n + 2m^2
+ m & &\text{for $\ba \in \Psi'_{2n,m}$,} \\
0 &\le a_0, & 0 &\le a_1, & a_i &\le a_{i+2}-1, & \sum a_i &= n + m^2 
& &\text{for $\ba \in \Phi'_{2n,m}$}
\end{align*}
There are embeddings $S: \Psi'_{2n,m} \to \Psi'_{2n,m+1}$, $S:
\Phi'_{2n,m} \to \Phi'_{2n,m+1}$ given by
\begin{align*}
S(\ba) &= (0,1,a_0+2,\ldots,a_{2m}+2) & &\text{for $\ba \in
  \Psi'_{2n,m}$,} \\
S(\ba) &= (0,0,a_0+1,\ldots,a_{2m}+1) & &\text{for $\ba \in
  \Phi'_{2n,m}$.} 
\end{align*}
These maps are called \emph{shift operations}.  Put $\Psi'_{2n} = \lim_{m\to \infty} \Psi'_{2n,m}$ and $\Phi'_{2n} =
\lim_{m\to\infty} \Phi'_{2n,m}$.  These sets are finite; we fix an $m$
large enough that $\Psi'_{2n} \simeq \Psi'_{2n,m}$ and $\Phi'_{2n} \simeq
\Phi'_{2n,m}$.  Elements of $\Phi'_{2n}$ (resp.~$\Psi'_{2n}$) are called 
\emph{symbols} (resp.~\emph{$u$-symbols}).

The bijection $i: \Psi'_{2n} \to \Phi'_{2n}$ is given by
\[
i(a_0, \ldots, a_{2m}) = (a_0, a_1 -1, a_2 -1, \ldots, a_{2m-1}-m, a_{2m}-m).
\]

\begin{defn}
  A \emph{block} (for the distinguished $u$-symbol $\ba$) is a subset
  $[k,l]$ of $[0,2m]$ that satisfies one of the following conditions:
\begin{itemize}
\item $[k,l]$ is a union of consecutive parts $P_1, \ldots, P_r$, where
$P_1$ and $P_r$ are ladders, $P_r$ is the unique part with even top, and
$P_1$ is the unique part with odd bottom. 
\item $k = 0$, and $[0,l]$ is the union of consecutive parts $P_2, \ldots,
P_r$, where $P_r$ is a ladder and the only part with even top and where
all parts have even bottom. 
\end{itemize}
In either case, the \emph{top} of the block is the top of $P_r$, and
its \emph{bottom} is the bottom of $P_1$, provided the latter is
defined.  Blocks satisfying the second condition above are said to
have no bottom.  Note that there is always exactly one block of this
type.
\end{defn}

We remark that it is permitted that $r = 1$, in which case the block
consists of a single ladder.  A part belongs to some block if and only if
it is not a staircase with odd bottom (and hence even top). 

\begin{lem}\label{lem:isol}
For any $\ba \in \Psi'_{2n}$, there are an odd number of isolated points. 
Suppose $i(\ba) = (a'_0, \ldots, a'_{2m})$, and let $k_0, \ldots, k_{2f}$
be the isolated points, numbered such that
\[
a'_{k_0} < a'_{k_1} < \cdots < a'_{k_{2f}}.
\]
If $i(\ba)$ is special, then $k_t \equiv t \pmod{2}$.
\end{lem}
\begin{proof}
  It is obvious from the definition of symbol that the even and odd
  indexed entries are strictly increasing, so that an entry can appear
  at most twice.  Hence, the number of isolated points is odd; it is
  the length of the symbol minus twice the number of duplicated
  entries.

Now suppose that $i(\ba)$ is special.  Every entry before $a'_{k_0}$
appears twice, so $k_0$ is even.  Similarly, every entry between
$a'_{k_t}$ and $a'_{k_{t+1}}$ is repeated, so $k_{t+1}\equiv k_t + 1
\equiv t+1 \pmod{2}$ by induction.
\end{proof}

We will adhere to the convention introduced in this lemma for naming
the isolated points.  In particular, this means that if $i(\ba)$ is not special, it
is not necessarily the case that $k_t < k_{t+1}$.

If $i(\ba)$ is a special symbol, we identify $V(i(\ba))$ with
the space $V_\bc$ of Section~\ref{sect:Vsp} (where $\bc$ is the
two-sided cell corresponding to $i(\ba)$) by 
\[
e_i = \{k_{i-1}, k_i\}, \qquad i = 1, \ldots, 2f;
\qquad\qquad
e_0 = \{k_1, k_2, \ldots, k_{2f}\}.
\]

\begin{lem}\label{lem:details}
Let $\ba = (a_0, \ldots, a_{2m})$ be a distinguished $u$-symbol, and let
$\ba' = i(\ba) = (a'_0, \ldots, a'_{2m})$ be the corresponding symbol.  Let
$k$ and $l$ be two integers such that $0 \le k < l \le 2m$. 
\begin{enumerate}
\item If $k$ and $l$ belong to distinct parts, then $a'_k <
  a'_l$.\label{pt:dist} 
\item If $[k,l]$ is a staircase with even bottom, then $k+1$ and $l-1$
  are the only isolated points in $[k,l]$.  A staircase with odd
  bottom contains no isolated points.\label{pt:staircase}
\item If $[k,l]$ is a ladder, then $k$ is an isolated point if it is
  odd, and $l$ is an isolated point if it is even.  There are no
  other isolated points in $[k,l]$.\label{pt:ladder}
\item The number of isolated points in any part is congruent to the 
  length of the part modulo $2$.\label{pt:count} 
\item For each $t$, either $k_t < k_{t+1}$, or $k_t = k_{t+1} +1$ and
  $[k_{t+1}, k_t]$ is a staircase with even bottom.\label{pt:order}
\item Displaced isolated points occur in pairs $\{k_t,k_{t+1}\} =
  e_{t+1}$ with $t$ even, with one such pair for each staircase with
  even bottom.\label{pt:disp}
\end{enumerate}
\end{lem}
\begin{proof}
\begin{parts}
\item If $k \equiv l \pmod{2}$ but $k \ne l$, it follows from the
  definition of $\Phi'_{2n}$ that $a'_k < a'_l$.  Indeed, the fact
  that $a'_i \le a'_{i+2} -1$ implies that $a'_i \le a'_{i+2j} -j$.  We
  need only treat the case where $k \not\equiv l \pmod{2}$.

Assume now that $k < l$, $k$ is even, and $l$ is odd.  Since $a_{l-1}
\le a_l$, $a'_{l-1} = a_{l-1} - (l-1)/2$, and $a'_l = a_l - (l+1)/2$,
it follows that $a'_l \ge a'_{l-1} - 1$.  Therefore, assuming $k <
l-3$, we have
\[
a'_k \le a'_{l-1} - (l-1-k)/2 < a'_{l-1} - 1 \le a'_l.
\]
If $k = l-1$, the fact that $k$ and $l$ are in different
parts implies that $a_{l-1} \le a_l -2$, whence $a'_l \ge a'_{l-1} +
1$.  Thus, $a'_k < a'_l$.  

Finally, suppose $k=l-3$.  If $l-1$ and $l$ are in different parts,
then the previous case gives $a'_l \ge a'_{l-1} +
1 \ge a'_{l-3}+2$.  Similar reasoning applies if $l-3$ and $l-2$ are
in different parts.  Otherwise, we have $l-3$ and $l-2$ in one part
and $l-1$ and $l$ in another.  This means that $(a_{l-3},\dots,a_l)$
has the form $(a,a,b,b)$, $(a-1,a,b,b)$, or  $(a,a,b,b+1)$ with $b-a>2$
or $(a-1,a,b,b+1)$ with $b-a>1$.  In each case, $a_l> a_{l-3}+2$, so
that $a'_l=a_l-(l+1)/2 > a_{l-3}-(l-3)/2 = a'_{l-3}$.

The reasoning is similar if $k$ is odd and $l$ is even.

\item Suppose $(a_k, a_{k+1}, \ldots, a_l) = (a,a,a+2,a+2, \ldots,
a+l-k-1)$.  For $i \in [k,l]$, it is easy to verify that:
\[
a'_i = 
\begin{cases}
a - k + i/2 & \text{if $k$ is even and $i$ is even,} \\
a - k + (i-3)/2& \text{if $k$ is even and $i$ is odd,} \\
a - k + (i-2)/2 & \text{if $k$ is odd and $i$ is even,} \\
a - k + (i-1)/2 & \text{if $k$ is odd and $i$ is odd.}
\end{cases}
\]
If $k$ is odd, then $a'_i = a'_{i+1}$ for $i = k, k+2, \ldots,
l-1$, so there are no isolated points.  On the other hand, if $k$ is
even, then $a'_i = a'_{i+3}$ for $i = k, k+2, \ldots, l-3$, so
none of the integers
\[
k, k+2, \ldots, l-3;\quad k+3, k+5, \ldots, l
\]
are isolated.  However, $a'_{k+1}$ and $a'_{l-1}$ are not duplicated.
In view of part~(\ref{pt:dist}), we see that $k+1$ and $l-1$ are
isolated.

\item Suppose $(a_k, a_{k+1}, \ldots, a_l) = (a, a+1, \ldots,
  a+l-k)$.  For $i \in [k,l]$, we have
\[
a'_i =
\begin{cases}
a - k + i/2 & \text{if $i$ is even,} \\
a - k + (i-1)/2 & \text{if $i$ is odd.}
\end{cases}
\]
Thus, $a'_i = a'_{i+1}$ if $i$ is even and  $i,i+1 \in [k,l]$.  It follows
that $k$ is isolated if and only if it is odd, that $l$ is isolated if
and only if it is even, and that no other points can be isolated.

\item From the preceding parts of the lemma, we see that each
staircase and each ladder of even length contains $0$ or $2$
isolated points, whereas each ladder of odd length contains
exactly one isolated point.

\item From part~(\ref{pt:dist}) we see that the inequality $k_t <
k_{t+1}$ can be violated only if $k_t$ and $k_{t+1}$ belong to the
same part.  Examining the formulas above for isolated points in
ladders and staircases yields the result.

\item From parts~(\ref{pt:ladder}) and~(\ref{pt:count}), it follows that
  no isolated points in ladders can be displaced, while
  parts~(\ref{pt:staircase}) and~(\ref{pt:count}) together imply that all
  isolated points in staircases are displaced.  Finally, if $[k,l]$ is
  a staircase with even bottom, then $k_t=k+1$, so that $k_t$ is odd.
  Since a special symbol has $k_t\equiv t \pmod{2}$, we see that $t$
  is even.
\end{parts}
\end{proof}

The following proposition is now an immediate consequence of the
definition of blocks and parts~(\ref{pt:staircase}) and~(\ref{pt:ladder})
of the preceding lemma. 

\begin{prop}\label{prop:isolC}
Let $\ba$ be a distinguished $u$-symbol. The top and bottom of each
block of $\ba$ are nondisplaced isolated points.  In addition to these,
each block contains an even number (possibly zero) of additional isolated
points, all of which are displaced.  There are no isolated points that do
not belong to any block. \qed 
\end{prop}

\begin{prop}
  For any block $b$ of a distinguished $u$-symbol $\ba$, the element
  $m_b \in \cG_\bc \times \hat \cG_\bc$ is of the form $(x,1)$ for
  some $x \in \cG_\bc$.
\end{prop}
\begin{proof}
It follows immediately from Lemma~\ref{lem:details}(\ref{pt:disp}) that
$\tilde\kappa(\ba)$ is a sum of $e_i$'s with $i$ odd.  Thus, $\tilde\kappa(\ba) \in
\cG_\bc \times \{1\} \subset M(\cG_\bc)$. 
\end{proof}

If $C$ is the unipotent class corresponding to $\ba$, then we set
$t(C)$ equal to the $x$ appearing in the above proposition.

\begin{prop}
If $C_1$ and $C_2$ are two unipotent classes in the same special piece,
then $t(C_1) \le t(C_2)$ in the partial order on $\cG_\bc$ if and only if 
$C_1$ is in the closure of $C_2$.
\end{prop}
\begin{proof}
Let $\ba$ be a distinguished $u$-symbol corresponding to a unipotent class
$C$, and suppose that $k_{t-1}$ and $k_t$ are consecutive nondisplaced
isolated points of $\ba$ with $t$ odd.  Thus, if we write $\tilde\kappa(\ba)
= (x,1)$, then the basis element $e_t \in \cG_\bc$ does not occur in $x$. 
Now, assume that the twisted symbol $i(\ba)^{e_t}$ is defined, and let
$\ba'$ be the $u$-symbol satisfying $i(\ba') = i(\ba)^{e_t}$.  Furthermore,
let us assume that $\ba'$ is also distinguished, corresponding to the
unipotent class $C'$.  Evidently, we have
$\tilde\kappa(\ba') = (x + e_t, 1)$ and $x+e_t < x$.  The first step will be
to prove that in this context, $C'$ lies in the closure of $C$.  We
introduce the following notation: 
\begin{align*}
\ba &= (a_0, \ldots, a_{2m}) &
i(\ba) &= \bb = (b_0, \ldots, b_{2m}) \\
\ba' &= (a'_0, \ldots, a'_{2m}) &
i(\ba') = \bb^{e_t} &= \bb' = (b'_0, \ldots, b'_{2m})
\end{align*}

Now, from the description of isolated points in
Lemma~\ref{lem:details} and Proposition~\ref{prop:isolC}, we know that
$k_{t-1}$ must be the top of a block and $k_t$ the bottom of the next
block.  Therefore, $[k_{t-1}+1, k_t-1]$ is a (possibly empty)
union of staircases with odd bottom.  In particular, $a_i = a_{i+1}$
for all odd $i \in [k_{t-1}+1, k_t - 1]$, and the
formula for the bijection $i$ gives $b_i = b_{i+1}$ for all odd $i
\in [k_{t-1}+1, k_t-1]$ as well.  (On the other hand, we note for
later reference that $a_i \le a_{i+1}-2$ for all even $i \in [k_{t-1},
k_t]$.)  Using these observations, one can show that the twisted
symbol is described by
\[
b'_i = \begin{cases}
b_i   & \text{if $i \notin [k_{t-1}, k_t]$,} \\
b_{i+1} & \text{if $i \in [k_{t-1}, k_t]$ and $i$ is even,} \\
b_{i-1} & \text{if $i \in [k_{t-1}, k_t]$ and $i$ is odd.}
\end{cases}
\]
(It suffices to verify that the sequence $(b_0, \ldots, b_{2m})$ is a valid
symbol and that it satisfies the parity conditions in the definition of
twisting.) 

Again using the formula for $i$, we find that
\[
a'_i = \begin{cases}
a_i     & \text{if $i \notin [k_{t-1}, k_t]$,} \\
a_{i+1}-1 & \text{if $i \in [k_{t-1}, k_t]$ and $i$ is even,} \\
a_{i-1}+1 & \text{if $i \in [k_{t-1}, k_t]$ and $i$ is odd.}
\end{cases}
\]
Since $a_i \le a_{i+1}-2$ for all even $i \in [k_{t-1},k_t]$, we
observe in particular that if $i \in [k_{t-1},k_t]$ is even, then 
\[
a'_i > a_i,
\qquad
a'_{i+1} < a_{i+1},
\qquad\text{and}\quad
a'_i+a'_{i+1} = a_i + a_{i+1}.
\]
Now, for any $u$-symbol $\bc = (c_0, \ldots, c_{2m})$, we define the sums
\[
\sigma_j(\bc) = \sum_{i=j}^{2m} c_i.
\]
From the above formulas, it is easy to see that
\[
\begin{aligned}
\sigma_j(\ba') &< \sigma_j(\ba) && \text{if $j \in [k_{t-1},k_t]$ is
  odd, and} \\
\sigma_j(\ba') &= \sigma_j(\ba) && \text{for all other $j \in [0,2m]$.}
\end{aligned}
\]
It was shown in~\cite[Lemme~3.2]{AA} that $C'$ is contained in the closure
of $C$ if and only if $\sigma_j(\ba') \le \sigma_j(\ba)$ for all $j$.  So
$C'$ is indeed contained in the closure of $C$. 

More generally, if $C_1$ and $C_2$ are two unipotent classes in the
same special piece with the property that $t(C_1) < t(C_2)$, then
$t(C_1)$ contains various terms $e_i$ (with $i$ odd) not appearing in
$t(C_2)$.  Adding these terms to $t(C_2)$ one at a time and iterating
the above argument shows that $C_1$ is contained in the closure of
$C_2$.

It remains to prove the opposite implication.  We will show that if
$t(C_1)$ and $t(C_2)$ are incomparable, then so are $C_1$ and $C_2$.
Let $\ba_1$ and $\ba_2$ be the corresponding distinguished
$u$-symbols, and let $\ba_0$ be the $u$-symbol corresponding to the
special class in their special piece.  Let $k_0, \ldots, k_{2f}$ be
the isolated points of $\ba_0$.  Since $t(C_1)$ is not contained in
$t(C_2)$, there exists a pair of consecutive isolated points
$\{k_{s-1}, k_s\}$ with $s$ odd, such that the corresponding isolated
points of $\ba_1$ are displaced, but those of $\ba_2$ are not.  In
other words, the term $e_s$ appears in $t(C_1)$, but not in $t(C_2)$.
On the other hand, since $t(C_2)\nsubseteq t(C_1)$, there exists
another such pair $\{k_{t-1}, k_t\}$ with $t$ odd, which is displaced
in $\ba_2$, but not in $\ba_1$.  Choose an odd number $j_1$ between
$k_{s-1}$ and  $k_s$ and another odd number $j_2$ between $k_{t-1}$
and $k_t$.  The
above calculations of $\sigma_j$ show that
\begin{align*}
\sigma_{j_1}(\ba_1) < &\sigma_{j_1}(\ba_0) = \sigma_{j_1}(\ba_2),
\\
\sigma_{j_2}(\ba_1) = &\sigma_{j_2}(\ba_0) > \sigma_{j_2}(\ba_2).
\end{align*}
Again invoking \cite[Lemme~3.2]{AA}, we see that this pair of inequalities
implies that neither $C_1$ nor $C_2$ may be contained in the closure of the
other. 
\end{proof}

\begin{prop}
Let $\ba$ be a distinguished $u$-symbol, corresponding to the unipotent
class $C$.  The elements $\pi(m_b)$, as $b$ runs over the blocks of $\ba$, are
  the components of $t(C)$, and the elements $\pi(\lambda_b)$ are the
  images in $V_\bc$ of the blocks of $t(C)$.
\end{prop}
\begin{proof}
  Let $b$ a block of $\ba$.  First, assume that $b$ has a bottom.  Let
  $k_{t-1}, k_t, \ldots, k_{t+2j}$ be the isolated points of $b$,
  where $k_{t-1}$ is the bottom of $b$ and $k_{t+2j}$ is its top.
  Since $k_{t-1}$ is odd and nondisplaced, $t$ is even.  Moreover,
  $k_{t-2}$ cannot be displaced, because it must be the top of the
  preceding block, while $k_{t+2j+1}$ is either undefined or the
  bottom of the next block.  This implies that the pairs $\{k_{t-2},
  k_{t-1}\}$ and $\{k_{t+2j}, k_{t+2j+1}\}$ contain no displaced
  isolated points, so the terms $e_{t-1}$ and $e_{t+2j+1}$ (provided
  the latter is defined) do not appear in $\tilde\kappa(\ba)$.
  Accordingly, $e_t + e_{t+2} + \cdots + e_{t+2j}$ is a block of
  $t(C)$ with associated component $e_{t+1} + e_{t+3} + \cdots + e_{
t+2j-1}$.  Projecting to $V_\bc$, we obtain
\[
\begin{aligned}
\pi(\lambda_b) &= \pi(\{k_{t-1}, \ldots, k_{t+2j}\}) = e_t + e_{t+2} + \cdots +
e_{t+2j}, \\
\pi(m_b) &= \pi(\{k_{t}, \ldots, k_{t+2j-1}\}) = e_{t+1} + e_{t+3} + \cdots + e_{
t+2j-1},
\end{aligned}
\]
where we abuse notation slightly by viewing the $e_i$'s as elements of
$V_\bc$.

If $b$ has no bottom, then the isolated points of $b$ are of the form
$k_0, k_1, \ldots, k_{2j}$.  Just as before, $e_0 + e_2 + \cdots +
e_{2j}$ is a block of $t(C)$ with associated component $e_{1} + e_{3}
+ \cdots + e_{2j-1}$ while the above formula for $\pi(m_b)$ is still
valid, with $t = 0$.  Recall,however,  that the definition of
$\lambda_b$ is different when $b$ contains an odd number of isolated
points: if $\ba$ contains $2f+1$ isolated points in all, then
\begin{equation*}
\pi(\lambda_b) = \pi(\{k_{2j+1}, k_{2j+2}, \ldots, k_{2f}\}) = e_{2j+2} + e_{2j+4} +
\cdots + e_{2f}.
\end{equation*}
Nevertheless, this is still the image in $V_\bc$ of a block, namely $e_0 + e_2 + \cdots + e_{2j}$.
\end{proof}

\begin{prop} The set $\cT(\ba)$ of elements $b\in\Hb(\ba)$ for which
  $\ba^b$ is defined is the codimension one hyperplane consisting of
  all subsets not containing the unique bottomless block.  Moreover,
  for all such $b$, $i(\ba^b) = (i(\ba))^{\lambda_b}$.
\end{prop}
\begin{proof}
It suffices to prove that $\ba^b$ is defined and that the formula for
$i(\ba^b)$ holds for the case where $b=[k,l]$ is a single block with a bottom.
We introduce the following notation:
\begin{align*}
\ba &= (a_0, \ldots, a_{2m}) &
i(\ba) &= (b_0, \ldots, b_{2m}) \\
\ba^b &= (a'_0, \ldots, a'_{2m}) &
i(\ba)^{\lambda_b} &= (b'_0, \ldots, b'_{2m}) \\
i(\ba^b) &= (c_0, \ldots, c_{2m})
\end{align*}
Now, from the definition of block, it follows that $a_{k-1}+1 < a_i <
a_{l+1}-1$ for all $i \in [k,l]$ (as long as $a_{l+1}$
is defined).  From this observation, one can deduce
that 
\[
a'_i = \begin{cases}
a_i & \text{if $i \notin [k,l]$,} \\
a_{i-1} & \text{if $i \in [k,l]$ and $i$ is even,} \\
a_{i+1} & \text{if $i \in [k,l]$ and $i$ is odd.} 
\end{cases}
\]
Indeed, it suffices to check that the sequence defined by this formula
is a valid $u$-symbol and that it satisfies the parity conditions in
the definition of twisting.  Both of these are easy to verify.

Next, it follows from the calculations in the proofs of
parts~(\ref{pt:staircase}) and~(\ref{pt:ladder}) of Lemma~\ref{lem:details} that
\[
\{ b_i \mid \text{$i \in [k,l]$ is odd and not isolated} \} =
\{ b_i \mid \text{$i \in [k,l]$ is even and not isolated} \}.
\]
  This observation lets us conclude that $i(\ba)^{\lambda_b}$ is
  defined and 
\[
b'_i = \begin{cases}
b_i & \text{if $i \notin [k,l]$,} \\
b_{i-1} & \text{if $i \in [k,l]$ and $i$ is even,} \\
b_{i+1} & \text{if $i \in [k,l]$ and $i$ is odd.} 
\end{cases}
\]

Finally, from the formula for $i$, we have
\[
b_i = 
\begin{cases}
a_i - i/2 \\
a_i - (i+1)/2
\end{cases}
\qquad\text{and}\qquad
c_i = 
\begin{cases}
a'_i - i/2     & \text{if $i$ is even,} \\
a'_i - (i+1)/2 & \text{if $i$ is odd.}
\end{cases}
\]
It follows that
\begin{align*}
c_i &=
\begin{cases}
a_i - i/2       & \text{if $i \notin [k,l]$ and $i$ is even,} \\
a_i - (i+1)/2   & \text{if $i \notin [k,l]$ and $i$ is odd,} \\
a_{i-1} - i/2   & \text{if $i \in [k,l]$ and $i$ is even,} \\
a_{i+1}-(i+1)/2 & \text{if $i \in [k,l]$ and $i$ is odd}
\end{cases}
\\
&\qquad\qquad\qquad\qquad
= 
\begin{cases}
b_i & \text{if $i \notin [k,l]$,} \\
b_{i-1} & \text{if $i \in [k,l]$ and $i$ is even,} \\
b_{i+1} & \text{if $i \in [k,l]$ and $i$ is odd} 
\end{cases}
\quad
= b'_i.
\end{align*}
This is the desired equality.
\end{proof}

\begin{prop}
The image of $p|_{\Hb(\ba)}: \Hb(\ba) \to \Irr(A(C))$ is precisely
$\Irr(\Ab(C))$.
\end{prop}
\begin{proof}
Recall that $\Ab(C)$ is defined to be the smallest quotient of $A(C)$ such
that all local systems on $C$ that are assigned by the Springer
correspondence to representations in the same two-sided cell as
$\chi_{C,1}$ actually come from representations of $\Ab(C)$.  Equivalently,
$\Irr(\Ab(C))$ can be characterized as the smallest subgroup of
$\Irr(A(C))$ that enjoys the following property:
\begin{quote}
For any $u$-symbol $\ba'$ that is similar to $\ba$, $\nu(\ba') \in
\Irr(\Ab(C))$ if $i(\ba')$ is congruent to $i(\ba)$.
\end{quote}
We first show that the image of $p|_{\Hb(\ba)}$ lies in $\Irr(A(C))$.
Since $\cT(\ba)$ is a transversal of the kernel, it suffices to show
that $p(\mu)\in\Irr(A(C))$ for each $\mu\in\cT(\ba)$.  For such a
$\mu$, the $u$-symbol $\ba^\mu$ is defined and is similar to $\ba$.
Moreover, the preceding proposition tells us that $i(\ba^\mu) =
i(\ba)^{\lambda_\mu}$.  Thus, $i(\ba^\mu)$ and $i(\ba)$ are congruent,
so $p(\mu)=\nu(\ba^\mu)$ lies in $\Irr(\Ab(C))$.

We use a dimension argument to prove that the image of $p|_{\Hb(\ba)}$
is all of $\Irr(\Ab(C))$.  It is a consequence of
Proposition~\ref{prop:kernel} that the image is an $\F_2$-vector space
whose dimension is given by any of the following:
\begin{gather*}
\text{(number of blocks in $\ba$)} - 1
\quad = \quad
\text{(number of blocks with a bottom)}\\
\ \quad = \quad
\text{(number of ladders with odd bottom)}
\end{gather*}
We now show that this is also the dimension of $\Ab(C)$.  Recall that
unipotent classes in type $C$ are in one-to-one correspondence with
partitions of $2n$ in which odd parts occur with even multiplicity.  Let
$\lambda$ be the partition corresponding to the same unipotent class
to which
$\ba$ corresponds.  Write $\lambda$ as 
\[
\lambda = (0 \le \lambda_0 \le \cdots \le \lambda_{2m})
\]
by increasing $m$ or adding zero parts as necessary.  In the course of the
proof of~\cite[Lemme~3.2]{AA}, a direct formula for $\ba$ in terms of
$\lambda$ is obtained:
\begin{equation}\label{eqn:ba-lambda}
a_i =
\begin{cases}
\frac{1}{2}\lambda_i + i + 1 & \text{if $\lambda_i$ is even,} \\
\frac{1}{2}\lambda_i + i + \frac{1}{2} &
  \text{if $\lambda_i$ is odd and $|\{j > i \mid \text{$\lambda_j$
  odd}\}|$ is even,}\\
\frac{1}{2}\lambda_i + i + \frac{3}{2} &
  \text{if $\lambda_i$ is odd and $|\{j > i \mid \text{$\lambda_j$
  odd}\}|$ is odd.}
\end{cases}
\end{equation}
(In~\cite{AA}, the entries of $\lambda$ and of $\ba$ (there called $\mu$)
are indexed by $i$ with $1 \le i \le 2m+1$.  Since we index them here by
$i$ with $0 \le i \le 2m$, we must add $1$ to the formulas
given there for $\mu_i$.)  Now, consider a fragment of the partition:
\[
\cdots \le \lambda_{k-1} < \lambda_k = \cdots = \lambda_l < \lambda_{l+1}
\le \cdots.
\]
Let $p$ be the common value of $\lambda_k, \ldots, \lambda_l$.  If $p$ is
even, the above formula for $\ba$ makes it clear that $[k,l]$ is a ladder.
On the other hand, if $p$ is odd, the fact that every odd part has even
multiplicity means that $a_k$ must be calculated using the third case
in~\eqref{eqn:ba-lambda}, $a_{k+1}$ must be  calculated using the
second case, and so on.  From this, it is easy to see that $[k,l]$ must be a
staircase.  Now, define the height of the part $p$ to be
\[
|\{j \mid \lambda_j < p\}| + 1.
\]
(Here, we add $1$ to make our usage of ``height'' consistent
with~\cite{Som}; partitions in type $C$ are there assumed to
have an even number of parts, whereas we assume that they have
an odd number of parts.)  Evidently, the height of $p$ is simply $k+1$.

In particular, even parts of even height in $\lambda$ correspond to
ladders with odd bottom in $\ba$.  According to the last paragraph
of~\cite[Section~5]{Som}, the dimension of $\Ab(C)$ is precisely the number of
even parts of even height in $\lambda$.  Thus, the image of $p|_{\Hb(\ba)}$
must be the full character group $\Irr(\Ab(C))$.
\end{proof}

\subsection{Type $B$}

We let $\Psi_{2n+1,m}$ (resp.~$\Phi_{2n+1,m}$) denote the set of
sequences $\ba = (a_0, a_1, \ldots, a_{2m})$ satisfying the following
conditions:
\begin{align*}
0 &\le a_0, & 0 &\le a_1, & a_i &\le a_{i+2}-2, & \sum a_i &= n + 2m^2
 & &\text{for $\ba \in \Psi_{2n+1,m}$,} \\
0 &\le a_0, & 0 &\le a_1, & a_i &\le a_{i+2}-1, & \sum a_i &= n + m^2 
& &\text{for $\ba \in \Phi_{2n+1,m}$}.
\end{align*}
The shift operations $S: \Psi_{2n+1,m} \to \Psi_{2n+1,m+1}$, $S:
\Phi_{2n+1,m} \to \Phi_{2n+1,m+1}$ are given by
\begin{align*}
S(\ba) &= (0,0,a_0+2,\ldots,a_{2m}+2) & &\text{for $\ba \in
  \Psi_{2n+1,m}$,} \\
S(\ba) &= (0,0,a_0+1,\ldots,a_{2m}+1) & &\text{for $\ba \in
  \Phi_{2n+1,m}$.} 
\end{align*}
As before, we put $\Psi_{2n+1} = \lim_{m\to \infty} \Psi_{2n+1,m}$ and
$\Phi_{2n+1} = \lim_{m\to\infty} \Phi_{2n+1,m}$.  These sets are again
finite; we fix an $m$ large enough that $\Psi_{2n+1} \simeq \Psi_{2n+1,m}$
and $\Phi_{2n+1} \simeq \Phi_{2n+1,m}$.  Elements of $\Phi_{2n+1}$
(resp.~$\Psi_{2n+1}$) are called \emph{symbols}
(resp.~\emph{$u$-symbols}). 

The bijection $i: \Psi_{2n+1} \to \Phi_{2n+1}$ is given by
\[
i(a_0, \ldots, a_{2m}) = (a_0, a_1, a_2 -1, a_3-1,\ldots, a_{2m}-m).
\]

\begin{defn}
A \emph{block} (for the distinguished $u$-symbol $\ba$) is a subset $[k,l]$ of $[0,2m]$ that satisfies one of the
following conditions: 
\begin{itemize}
\item $[k,l]$ is a union of consecutive parts $P_1, \ldots, P_r$, where
$P_1$ and $P_r$ are ladders, $P_1$ is the unique part with even bottom, and
$P_r$ is the unique part with odd top. 
\item $l = 2m$, and $[k, 2m]$ is the union of consecutive parts $P_1,
\ldots, P_{r-1}$, where $P_1$ is a ladder and the only part with even
bottom and where all parts have even top. 
\end{itemize}
In either case, the \emph{bottom} of the block is the bottom of $P_1$,
and its \emph{top} is the top of $P_r$, provided the latter is
defined.  Blocks satisfying the second condition above are said to
have no top.  There is always exactly one block of this type.
\end{defn}

Note that a part belongs to some block if and only if it is not a staircase
with even bottom (and hence odd top).

\begin{lem} The set $\cT(\ba)$ of elements $b\in\Hb(\ba)$ for which
  $\ba^b$ is defined is the codimension one hyperplane consisting of
  all subsets not containing the unique topless block.
\end{lem}

\begin{lem}\label{lem:isolB}
For any $\ba \in \Psi_{2n+1}$, there are an odd number of isolated points. 
Suppose $i(\ba) = (a'_0, \ldots, a'_{2m})$, and let $k_0, \ldots, k_{2f}$
be the isolated points, numbered such that
\[
a'_{k_0} < a'_{k_1} < \cdots < a'_{k_{2f}}.
\]
If $i(\ba)$ is special, then $k_t \equiv t \pmod{2}$.
\end{lem}

As we did in type $C$, we will adhere to the convention introduced in
this lemma for naming the isolated points.  If $i(\ba)$ is a special
symbol, we identify $V(i(\ba))$ with the space $V_\bc$ of
Section~\ref{sect:Vsp} (where $\bc$ is the
two-sided cell corresponding to $i(\ba)$) by

\[
e_i = \{k_i, k_{i+1}\}, \qquad i = 0, \ldots, 2f-1;
\qquad\qquad e_{2f} = \{k_0, k_1, \ldots, k_{2f-1}\}.
\]

\begin{lem}\label{lem:detailsB}
Let $\ba = (a_0, \ldots, a_{2m})$ be a distinguished $u$-symbol, and let
$\ba' = i(\ba) = (a'_0, \ldots, a'_{2m})$ be the corresponding symbol.  Let
$k$ and $l$ be two integers such that $0 \le k < l \le 2m$. 
\begin{enumerate}
\item If $k$ and $l$ belong to distinct parts, then $a'_k <
  a'_l$. 
\item If $[k,l]$ is a staircase with odd bottom, then $k+1$ and $l-1$
  are the only isolated points in $[k,l]$.  A staircase with even
  bottom contains no isolated points.
\item If $[k,l]$ is a ladder, then $k$ is an isolated point if it is
  even, and $l$ is an isolated point if it is odd.  There are no
  other isolated points in $[k,l]$.
\item The number of isolated points in any part is congruent to the 
  length of the part modulo $2$.
\item For each $t$, either $k_t < k_{t+1}$, or $k_t = k_{t+1} +1$ and
  $[k_{t+1}, k_t]$ is a staircase with odd bottom.
\item Displaced isolated points occur in pairs $\{k_t,k_{t+1}\} =
  e_t$ with $t$ odd, with one such pair for each staircase with
  odd bottom.
\end{enumerate}
\end{lem}

\subsection{Type $D$}

Finally, let $\Psi_{2n,m}$ (resp.~$\Phi_{2n,m}$) denote the set of
sequences $\ba = (a_0, a_1, \ldots, a_{2m+1})$ satisfying the following
conditions:
\begin{align*}
0 &\le a_0, & 0 &\le a_1, & a_i &\le a_{i+2}-2, & \sum a_i &= n + 2m^2
+ 2m & &\text{for $\ba \in \Psi_{2n,m}$,} \\
0 &\le a_0, & 0 &\le a_1, & a_i &\le a_{i+2}-1, & \sum a_i &= n + m^2 +m
& &\text{for $\ba \in \Phi_{2n,m}$},
\end{align*}
together with the additional condition:
\begin{equation}\label{eqn:D} \text{if } i \text{ is the smallest
    index such that } a_{2i}\ne a_{2i+1}\text{, then }a_{2i}<
  a_{2i+1}.\tag{$\star$}\end{equation}
(Note that such an $i$ need not exist; in the case of $u$-symbols,
there is no such $i$ if and only if $\ba$ is a union of  staircases, necessarily all with even bottom.)

The shift operations $S: \Psi_{2n,m} \to \Psi_{2n,m+1}$, $S:
\Phi_{2n,m} \to \Phi_{2n,m+1}$ are given by
\begin{align*}
S(\ba) &= (0,0,a_0+2,\ldots,a_{2m+1}+2) & &\text{for $\ba \in
  \Psi_{2n,m}$,} \\
S(\ba) &= (0,0,a_0+1,\ldots,a_{2m+1}+1) & &\text{for $\ba \in
  \Phi_{2n,m}$.} 
\end{align*}
Put $\Psi_{2n} = \lim_{m\to \infty} \Psi_{2n,m}$ and $\Phi_{2n} =
\lim_{m\to\infty} \Phi_{2n,m}$.  We fix an $m$ large enough that $\Psi_{2n}
\simeq \Psi_{2n,m}$ and $\Phi_{2n} \simeq \Phi_{2n,m}$.  Elements of
$\Phi_{2n}$ (resp.~$\Psi_{2n}$) are called \emph{symbols}
(resp.~\emph{$u$-symbols}).

The bijection $i: \Psi_{2n} \to \Phi_{2n}$ is given by
\[
i(a_0, \ldots, a_{2m+1}) = (a_0, a_1, a_2 -1, a_3-1,\ldots, a_{2m}-m, a_{2m+1}-m).
\]

\begin{defn}
A subset $[k,l]$ of $[0,2m]$ is called a \emph{block} (for the
distinguished $u$-symbol $\ba$) if it is a union of
consecutive parts $P_1, \ldots, P_r$, where $P_1$ and $P_r$ are ladders,
$P_r$ is the unique part with odd top, and $P_1$ is the unique part with
even bottom.  The \emph{top} of the block is the top of $P_r$, and its
\emph{bottom} is the bottom of $P_1$. 
\end{defn}

A part belongs to some block if and only if it is not a staircase with even
bottom (and hence odd top). 

In type $D$, it is possible for a $u$-symbol $\ba$ to have no blocks
and no isolated points; this occurs precisely when $\ba$ is a union of
staircases.  In this situation, $\Ab(C)$ and $\cG'_\bc$ both have only
one element, and the main theorem is trivial.  For the rest of this
section, we assume where necessary that $\ba$ has a block  (and
accordingly, an isolated point).

The sequence of integers $\ba^\mu$ is defined for each
$\mu\in\Hb(\ba)$, but it is not necessarily a $u$-symbol as it need
not satisfy condition \eqref{eqn:D}.

\begin{lem} For each $\mu\in\Hb(\ba)$, exactly one of $\ba^\mu$ and
  $\ba^{\mu^c}$ satisfies condition \eqref{eqn:D}.  Thus,  $\cT(\ba)$ is
a transversal  of the subgroup  $\{\varnothing,B(\ba)\}$.
\end{lem}

\begin{lem}\label{lem:isolD}
For any $\ba \in \Psi_{2n}$, there are an even number of isolated
points.  Suppose $i(\ba) = (a'_0, \ldots, a'_{2m+1})$, and let $k_0,
\ldots, k_{2f+1}$ be the isolated points, numbered such that
\[
a'_{k_0} < a'_{k_1} < \cdots < a'_{k_{2f+1}}.
\]
If $i(\ba)$ is special, then $k_t \equiv t \pmod{2}$.
\end{lem}

If $i(\ba)$ is a special symbol, we identify $V(i(\ba))$ with
the space $\tilde V_\bc$ of Section~\ref{sect:Vsp} (where $\bc$ is the
two-sided cell corresponding to $i(\ba)$) by 
\[
e_i = \{k_i, k_{i+1}\}, \qquad i = 0, \ldots, 2f.
\]
Note that in type $D$, the map $\pi:V(i(\ba))\to V_\bc$ is not an isomorphism.

\begin{lem}\label{lem:detailsD}
Let $\ba = (a_0, \ldots, a_{2m+1})$ be a distinguished $u$-symbol, and let
$\ba' = i(\ba) = (a'_0, \ldots, a'_{2m+1})$ be the corresponding symbol. 
Let $k$ and $l$ be two integers such that $0 \le k < l \le 2m+1$. 
\begin{enumerate}
\item If $k$ and $l$ belong to distinct parts, then $a'_k <
  a'_l$.
\item If $[k,l]$ is a staircase with odd bottom, then $k+1$ and $l-1$
  are the only isolated points in $[k,l]$.  A staircase with even
  bottom contains no isolated points.
\item If $[k,l]$ is a ladder, then $k$ is an isolated point if it is
  even, and $l$ is an isolated point if it is odd.  There are no
  other isolated points in $[k,l]$.
\item The number of isolated points in any part is congruent to the 
  length of the part modulo $2$.
\item For each $t$, either $k_t < k_{t+1}$, or $k_t = k_{t+1} +1$ and
  $[k_{t+1}, k_t]$ is a staircase with odd bottom.
\item Displaced isolated points occur in pairs $\{k_t,k_{t+1}\} =
  e_t$ with $t$ odd, with one such pair for each staircase with
  odd bottom.
\end{enumerate}
\end{lem}

\end{document}